\documentclass{elsart}
\usepackage{amsfonts}
\usepackage{amssymb}
\usepackage{graphicx}

\newcommand{\fin}{\hspace{\fill} \qed }

\begin{document}

\begin{verbatim}\end{verbatim}\vspace{2.5cm}

\begin{frontmatter}
\title{Characterizing [h,2,1] graphs by minimal forbidden induced subgraphs}
\author[laplata]{L. Alc\'on}
\ead{liliana@mate.unlp.edu.ar}
\author[laplata,coni]{M. Gutierrez}
\ead{marisa@mate.unlp.edu.ar}
\author[laplata,coni]{M. P. Mazzoleni}
\ead{pia@mate.unlp.edu.ar}

\address[laplata]{Departamento de Matem\'atica,
Universidad Nacional de La Plata,  CC 172, (1900) La Plata,
Argentina}
\address[coni]{CONICET}

\begin{abstract}
An undirected graph $G$ is called a VPT graph if it is the vertex
intersection graph of a family of paths in a tree. The class of
graphs which admit a VPT representation in a host tree with
maximum degree at most $h$ is denoted by $[h,2,1]$. The classes
$[h,2,1]$ are closed by taking induced subgraphs, therefore each
one can be characterized by a family of minimal forbidden induced
subgraphs. In this paper we associate the minimal forbidden
induced subgraphs for $[h,2,1]$ which are VPT with (color)
$h$-critical graphs. We describe how to obtain minimal forbidden
induced subgraphs from critical graphs, even more, we show that
the family of graphs obtained using our procedure is exactly the
family of VPT minimal forbidden induced subgraphs for $[h,2,1]$.
The members of this family together with the minimal forbidden
induced subgraphs for VPT \cite{B9,B11}, are the minimal forbidden
induced subgraphs for $[h,2,1]$, with $h\geq 3$. Notice that by
taking $h=3$ we obtain a characterization by minimal forbidden
induced subgraphs of the class VPT $\cap$ EPT=EPT $\cap$
Chordal$=[3,2,2]=[3,2,1]$ (see \cite{B4}).


\end{abstract}
\end{frontmatter}

\section{Introduction}

The \textbf{intersection graph} of a set family is a graph whose
vertices are the members of the family, and the adjacency between
them is defined by a non-empty intersection of the corresponding
sets. Classic examples are interval graphs and chordal graphs.\

An \textbf{interval graph} is the intersection graph of a family
of intervals of the real line, or, equivalently, the vertex
intersection graph of a family of subpaths of a path. A
\textbf{chordal graph} is a graph without chordless cycles of
length at least four. Gavril \cite{B1} proved that a graph is
chordal if and only if it is the  vertex intersection graph of a
family of subtrees of a tree. Both classes have  been widely
studied \cite{classes}.

In order to allow larger families of graphs to be represented by
subtrees, several graph classes are defined imposing conditions on
trees, subtrees and intersection sizes  \cite{Jami,B6}. Let $h$,
$s$ and $t$ be positive integers; an $(h,s,t)$-representation of a
graph $G$ consists in a host tree $T$ and a collection
$(T_v)_{v\in V(G)}$ of subtrees of $T$, such that (i) the maximum
degree of $T$ is at most $h$, (ii) every subtree $T_v$ has maximum
degree at most $s$, (iii) two vertices $v$ and $v'$ are adjacent
in $G$ if and only if the corresponding subtrees $T_v$ and
$T_{v'}$ have at least $t$ vertices in common in $T$. The class of
graphs that have an $(h,s,t)$-representation is denoted by
$[h,s,t]$. When there is no restriction on the maximum degree of
$T$ or on the maximum degree of the subtrees, we use $h=\infty$
and $s=\infty$ respectively. Therefore, $[\infty,\infty,1]$ is the
class of chordal graphs and $[2,2,1]$ is the class of  interval
graphs. The classes $[\infty,2,1]$ and $[\infty,2,2]$ are called
VPT and EPT  respectively in \cite{B4}; and UV and UE,
respectively in \cite{B5}.

 In \cite{gravril,nuevo}, it is shown that the problem of recognizing VPT graphs is polynomial time solvable.
Recently, in \cite{pia}, generalizing a result given in \cite{B4},
we have proved that the problem of deciding whether a given VPT
graph belongs to $[h,2,1]$ is NP-complete even when restricted to
the class VPT $\cap$ Split without dominated stable vertices. The
classes $[h,2,1]$, $h\geq 2$, are closed by taking induced
subgraphs, therefore each one can be characterized by a family of
minimal forbidden induced subgraphs. Such a family is known only
for $h=2$ \cite{B8} and there are some partial results for $h=3$
\cite{hugo}. In this paper we associate the VPT minimal forbidden
induced subgraphs for $[h,2,1]$ with (color) $h$-critical graphs.
We describe how to obtain minimal forbidden induced subgraphs from
critical graphs, even more, we show that the family of graphs
obtained using our procedure is exactly the family of VPT minimal
forbidden induced subgraphs for $[h,2,1]$. The members of this
family together with the minimal forbidden induced subgraphs for
VPT (see Figure \ref{fig:fiscpath})\cite{B9,B11}, are the minimal
forbidden induced subgraphs for $[h,2,1]$, with $h\geq 3$. Notice
that by taking $h=3$ we obtain a characterization by minimal
forbidden induced subgraphs of the class VPT$\cap$ EPT=EPT$\cap$
Chordal$=[3,2,2]=[3,2,1]$ \cite{B4}.

The paper is organized as follows: in Section $2$, we provide
basic definitions and basic results. In Section $3$, we give
necessary conditions for VPT minimal non $[h,2,1]$ graphs. In
Section $4$, we show a procedure to construct minimal non
$[h,2,1]$ graphs. In Section $5$, we characterize minimal non
$[h,2,1]$ graphs.

\section{Preliminaries}

\label{s:preliminares}  Throughout this paper,  graphs are
connected, finite and simple. The \textbf{vertex set}  and the
\textbf{edge set} of a graph $G$ are denoted by $V(G)$ and $E(G)$
respectively. The \textbf{open neighborhood} of a vertex $v$,
represented by $N_G(v)$, is the set of vertices adjacent to $v$.
The \textbf{closed neighborhood} $N_G[v]$ is $N_G(v)\cup \{v\}$.
The \textbf{degree} of $v$, denoted by $d_G(v)$, is the
cardinality of $N_G(v)$. For simplicity, when no confusion can
arise, we omit the subindex $G$ and write $N(v)$, $N[v]$ or
$d(v)$. Two vertices $x,y\in V(G)$ are called \textbf{true twins}
if $xy\in E(G)$ and $N(x)=N(y)$.\

A \textbf{complete set} is a subset of mutually adjacent vertices.
A \textbf{clique} is a maximal complete set. The family of cliques
of $G$ is denoted by $\mathcal{C}(G)$. A $\textbf{stable set}$ is
a subset of pairwise non-adjacent vertices.\

A graph $G$ is \textbf{k-colorable}  if its vertices can be
colored with at most $k$ colors in such a way that  no two
adjacent vertices share the same color. The \textbf{chromatic
number} of $G$, denoted by $\chi(G)$, is the smallest $k$ such
that $G$ is $k$-colorable. A vertex $v\in V(G)$ or an edge $e\in
E(G)$ is a $\textbf{critical element}$ of $G$ if
$\chi(G-v)<\chi(G)$ or $\chi(G-e)<\chi(G)$. A graph $G$ with
chromatic number $h$ is $\textbf{h-vertex critical}$ (resp.
$\textbf{h-edge critical}$) if each of its vertices (resp. edges)
is a critical element and it is $\textbf{h-critical}$ if both
hold.

A $\textbf{VPT representation}$ of $G$ is a pair $\langle
\mathcal{P},T \rangle$ where $\mathcal{P}$ is a family
$(P_v)_{v\in V(G)}$ of subpaths of a host tree $T$ satisfying that
two vertices $v$ and $v'$ of $G$ are adjacent if and only if $P_v$
and $P_{v'}$ have at least one vertex in common, in such case we
say that $P_v$ intersects $P_{v'}$. When the maximum degree of the
host tree is $h$ the VPT representation of $G$ is called an
$(h,2,1)$-representation of $G$. The class of graphs which admit
an $(h,2,1)$-representation is denoted by $\textbf{[h,2,1]}$.

Since a family of vertex paths in a tree satisfies the Helly
property \cite{berge}, if $C$ is a clique of $G$ then there exists
a vertex $q$ of $T$ such that $C=\{v\in V(G): q\in V(P_v) \}$. On
the other hand, if $q$ is any vertex of the host tree $T$, the set
$\{v\in V(G): q\in V(P_v) \}$, denoted by  $\mathbf{C_q}$, is a
complete set of $G$, but not necessarily a clique. In order to
avoid this drawback we introduce the notion of full representation
at $q$.

Let $\langle \mathcal{P},T\rangle$ be a VPT representation of $G$
and let $q$ be a vertex of degree $h$ of $T$. The connected
components of $T-q$ are called the \textbf{branches of T at q}. A
path is \textbf{contained} in a branch if all its vertices are
vertices of the branch. Notice that if $N_T(q)=\{q_1,q_2,..,q_h\}$
then $T$ has exactly $h$ branches at $q$. The branch containing
$q_i$ is denoted by $\mathbf{T_i}$. Two branches $T_i$ and $T_j$
are \textbf{linked} by a path $P_v\in \mathcal{P}$ if both
vertices $q_i$ and $q_j$ belong to $V(P_v)$.

\begin{defn} A VPT representation  $\langle \mathcal{P},T \rangle$ is
\textbf{full at a vertex q}  of $T$  if,  for every two branches
$T_i$ and $T_j$ of $T$ at $q$, there exist paths
   $ P_v,P_w,P_u\in \mathcal{P}$ such that: $(i)$ the branches $T_i$ and $T_j$ are linked
   by $P_v$; $(ii)$ $P_w$ is contained in $T_i$ and intersects $P_v$ in at least one vertex; and
   $(iii)$ $P_u$ is contained in $T_j$ and intersects $P_v$ in at least one
   vertex.\end{defn}

A clear consequence of the previous definition is that if $\langle
\mathcal{P},T \rangle$ is full at a vertex q of $T$, with
$d_T(q)=h\geq 3$, then $C_q$ is a clique of $G$.

The following theorem  shows that a VPT representation which is
not full at a vertex $q$ of $T$, with $d_{T}(q)=h\geq 4$, can be
modified to obtain a VPT representation without increasing the
maximum degree of the host tree; and, even more, decreasing the
degree of the vertex $q$.

\begin{thm}\cite{pia}
\label{t:repre} Let $\langle\mathcal{P},T\rangle$ be a VPT
representation of $G$. Assume there exists  a vertex $q\in V(T)$
with $d_{T}(q)=h\geq 4$ and two branches of $T$ at $q$ which are
linked by no path of $\mathcal{P}$.  Then there exists a VPT
representation $\langle \mathcal{P}',T'\rangle$ of $G$ with
$V(T')=V(T)\cup \{q'\}$, $q' \notin V(T)$, and

$$d_{T'}(x) =\left\{
\begin{array}{ll}
3,  \ \ \ &\mbox{if\ }x = q'\\
 \
h-1, \ \ \ &\mbox{if\ }x = q\\
 \
d_{T}(x), \ \ \ &\mbox{if\ }x \in V(T')\setminus \{q,q'\}.
\end{array}
\right.$$
\end{thm}

In what follows we give the definition of the branch graph which
can be used to describe intrinsic properties of representations of
VPT graphs.

\begin{defn}\cite{B4} \label{d:branch} Let $C\in \mathcal{C}(G)$. The \textbf{branch
graph}  of $G$ for the clique $C$, denoted by $\mathbf{B(G/C)}$,
is defined as follows: its vertices are the vertices of
$V(G)\setminus C$ which are adjacent to some vertex of $C$. Two
vertices $v$ and $w$ are adjacent in   $B(G/C)$ if and only if
\begin{enumerate}
\item $vw\notin E(G)$;

\item there exists a vertex $x\in C$ such that  $xv\in E(G)$ and
$xw\in E(G)$;

\item there exists a vertex $y\in C$ such that  $yv\in E(G)$ and
$yw \not\in E(G)$;

\item there exists a vertex $z\in C$ such that  $zv \not\in E(G)$
and $zw\in E(G)$.
\end{enumerate}
\end{defn}

It is clear that if $C\in \mathcal{C}(G)$ and $v\in V(G)-C$ then
$C\in \mathcal{C}(G-v)$. The following claim says what happens
with the branch graphs when we remove such vertices. Its proof is
trivial.

\begin{claim} \label{l:branch}
Let $C\in \mathcal{C}(G)$ and let $v\in V(G)-C$: $(i)$ If $v\notin
V(B(G/C))$ then $B(G-v/C)=B(G/C)$; (ii) if $v\in V(B(G/C))$ then
$B(G-v/C)=B(G/C)-v$.
\end{claim}

As will be seen in what follows, branch graphs of VPT graphs can
be used to describe intrinsic properties of representations.

\begin{lem}\cite{pia} \label{l:branch1} Let $C$ be a clique of a VPT graph $G$, $\langle
\mathcal{P}, T \rangle$ be a VPT representation of $G$ and $q$ be
a vertex of $T$ such that $C=C_q$. If $v$ is a vertex of $B(G/C)$
then $P_v$ is contained in some branch of $T$ at $q$. If two
vertices $v$ and $w$ are adjacent in $B(G/C)$ then  $P_v$ and
$P_{w}$ are not contained in a same branch of $T$ at $q$.
\end{lem}

In \cite{pia} we proved the following two results which show that
there is a relation between the VPT graphs that can be represented
in a tree with maximum degree at most $h$ and the chromatic number
of their branch graphs.

\begin{lem} \cite{pia}
\label{l:grado} Let $\langle \mathcal{P}, T \rangle$ be a VPT
representation of $G$. Let $C\in \mathcal{C}(G)$ and $q\in V(T)$
such that $C=C_{q}$. If $d_{T}(q)=h$, then $B(G/C)$ is
$h$-colorable.
\end{lem}

\begin{thm}\cite{pia}
\label{t:cromatic}
 Let $G\in$ VPT and  $h\geq 4$. The graph $G$ belongs to
 $[h,2,1]-[h-1,2,1]$ if and only if $Max_{C\in
\mathcal{C}(G)}(\chi(B(G/C)))=h$. The reciprocal implication is
also true for $h=3$.\end{thm}

\begin{defn} A clique $K$ of a graph $G$ is called \textbf{principal} if $$Max_{C\in
\mathcal{C}(G)}(\chi(B(G/C)))=\chi(B(G/K)).$$\end{defn}

A graph $G$ is \textbf{split} if $V(G)$ can be partitioned into a
stable set $S$ and a clique $K$. The pair $\mathbf{(S,K)}$ is the
\textbf{split partition} of $G$ and this partition is unique up to
isomorphisms. The vertices in $S$ are called \textbf{stable
vertices}, and $K$ is called the \textbf{central clique} of $G$.
We say that a vertex $s$ is \textbf{a dominated stable vertex} if
$s\in S$ and there exists $s'\in S$ such that $N(s)\subseteq
N(s')$. Notice that if $G$ is split then $\mathcal{C}(G)=\{K, N[s]
\mbox{ for } s\in S\}$. We will call \textbf{Split} to the class
of split graphs.

\begin{lem}\label{t:principal} Let $G\in$ VPT $\cap$ Split with split partition $(S,K)$. Then, $K$ is a
principal clique of $G$.\end{lem}

\begin{pf} Let $s\in S$, we know that $N[s]\in \mathcal{C}(G)$.
Observe that $V(B(G/N[s]))$ $=(K-N[s])\cup S'$, with $S'=\{x\in S:
N(x)\cap N(s)\neq \emptyset\}$. We claim that the vertices of
$K-N[s]$ are isolated in $B(G/N[s])$. Indeed, let $x\in K-N[s]$,
if $y\in K-N(s)$, then $xy\notin E(B(G/N[s]))$ because $xy\in
E(G)$ and, if $y\in S'$ then $xy\notin E(B(G/N[s]))$ because
$N(y)\subseteq N[x]$. Then, we are only interesting in the
subgraph of $B(G/N[s])$ induced by $S'$, and this is a subgraph of
$B(G/K)$. Thus, $\chi(B(G/N[s]))\leq \chi(B(G/K))$. Hence,
$Max_{C\in \mathcal{C}(G)}(\chi(B(G/C)))=\chi(B(G/K))$, that is,
$K$ is a principal clique of $G$.\fin\end{pf}

\section{Necessary conditions for VPT minimal non [h,2,1] graphs}
\label{s:caract}

In this Section we give some necessary conditions for VPT minimal
non [h,2,1] graphs, with $h\geq 3$; recall that:

\begin{defn} A \textbf{minimal non [h,2,1] graph} is a  minimal forbidden
induced subgraph for the class $[h,2,1]$, this means  any graph
$G$ such that $G \not\in [h,2,1]$ and $G-v \in [h,2,1]$ for every
vertex $v\in V(G)$.\end{defn}

\begin{thm}\label{t:contencion} Let $G\in$ VPT and let $h\geq 3$.
 If $G$ is a minimal non $[h,2,1]$ graph then $G\in [h+1,2,1]$.\end{thm}

\begin{pf} Let $C\in \mathcal{C}(G)$ and let $v\notin C$. We know
that $G-v\in [h,2,1]$ then, by Theorem \ref{t:cromatic},
$\chi(B(G-v/C))\leq h$. By Claim \ref{l:branch}, $\chi(B(G-v/C))=
\chi(B(G/C)-v)\geq$ $\chi(B(G/C))-1$. Thus, $\chi(B(G/C))-1\leq h$
and hence $\chi(B(G/C))\leq h+1$. Then, by Theorem
\ref{t:cromatic}, $G\in [h+1,2,1]$.\fin\end{pf}

\begin{thm} \label{l:adyacente1} Let $K$ be a principal clique of a VPT  minimal non
$[h,2,1]$ graph $G$, with $h\geq 3$.  Then: $(i)$ $V(B(G/K))=V(G)-
K$; $(ii)$ if $v\in V(G)-K$ then $|N(v)\cap K|>1$; $(iii)$
$B(G/K)$ is $(h+1)$-vertex critical; $(iv)$ if $s_1,s_2\in V(G)-K$
then $N(s_1)\cap K\neq N(s_2)\cap K$.\end{thm}

\begin{figure}[h!]
\footnotesize     \centering
\begin{tabular}{c}
\includegraphics[scale=0.35]{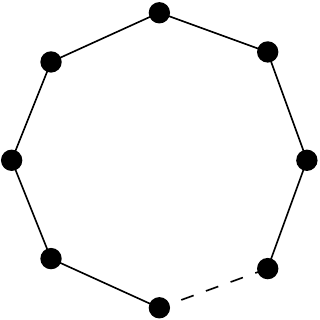} \\
$F_{0}(n)_{n\geq 4}$ \\
\end{tabular}

\begin{tabular}{ccccc}
\includegraphics[scale=0.35]{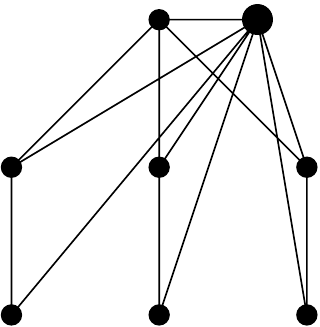} \ &\
\includegraphics[scale=0.35]{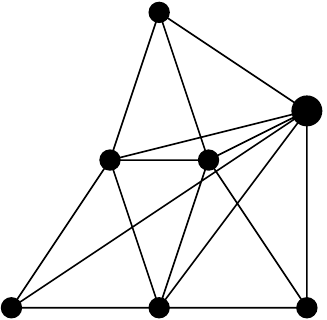} \ &\
\includegraphics[scale=0.35]{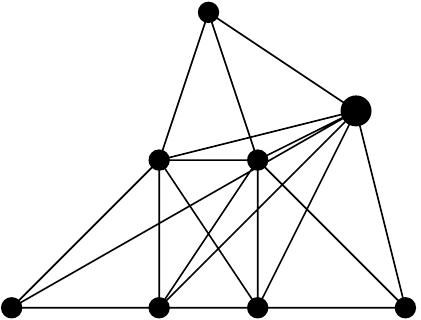} \ &\
\includegraphics[scale=0.35]{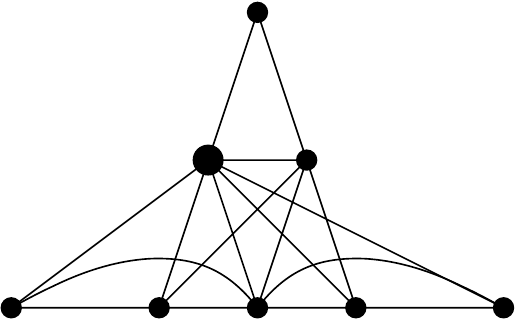} \ &\
\includegraphics[scale=0.35]{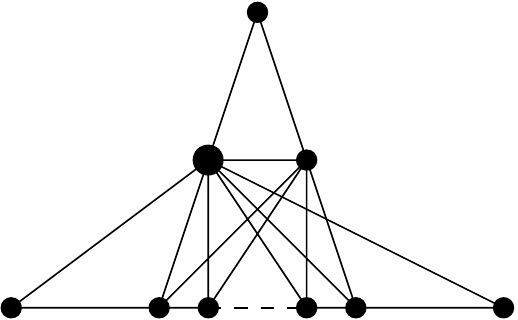} \\
$F_1$ & $F_2$ & $F_3$ & $F_4$ & $F_5(n)_{n\geq 7}$
\end{tabular}

\begin{tabular}{ccccc}
\includegraphics[scale=0.35]{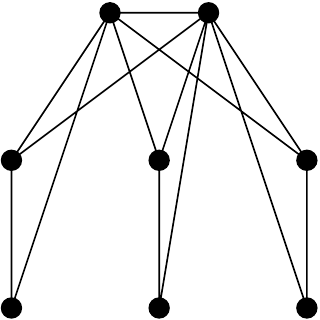} \ &\
\includegraphics[scale=0.35]{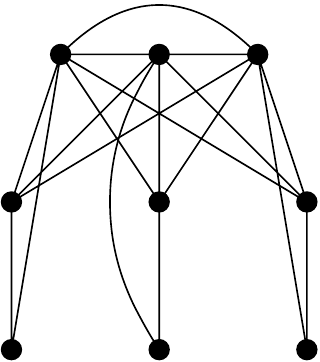} \ &\
\includegraphics[scale=0.35]{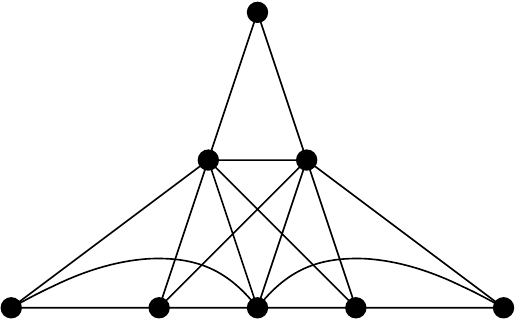} \ &\
\includegraphics[scale=0.35]{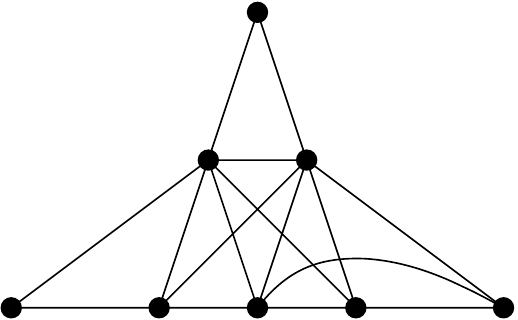}  \ &\
\includegraphics[scale=0.35]{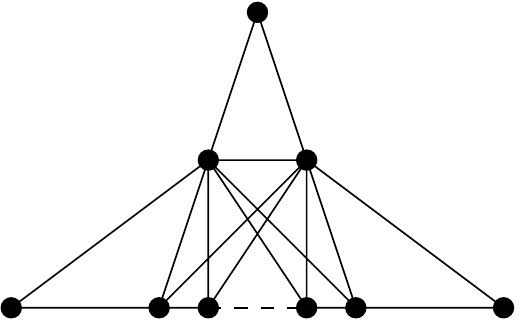} \\
$F_6$ & $F_7$ & $F_8$ & $F_9$ & $F_{10}(n)_{n\geq 8}$ \\
\end{tabular}

\begin{tabular}{cccccc}
\includegraphics[scale=0.35]{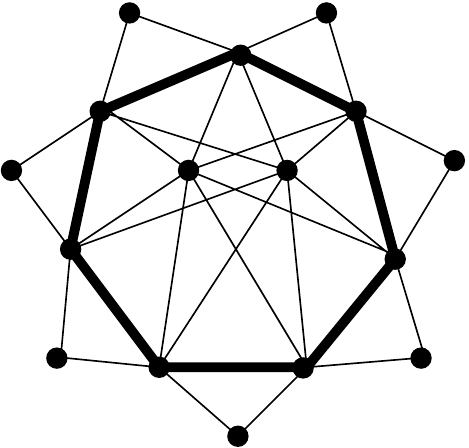} \ &\
\includegraphics[scale=0.35]{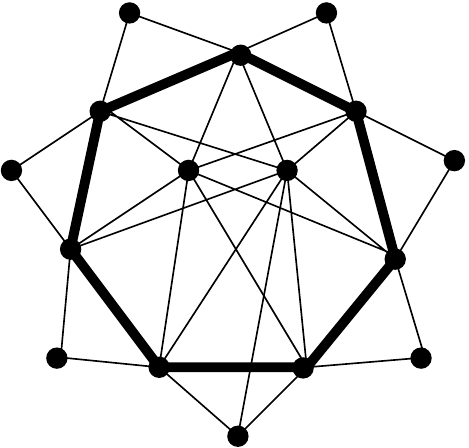} \ &\
\includegraphics[scale=0.35]{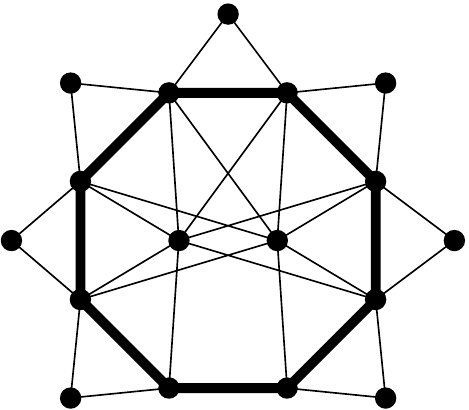}\ &\
\includegraphics[scale=0.35]{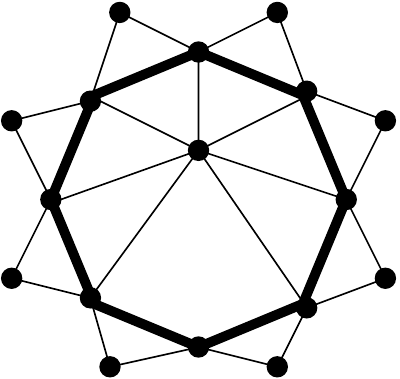}  \ &\
\includegraphics[scale=0.35]{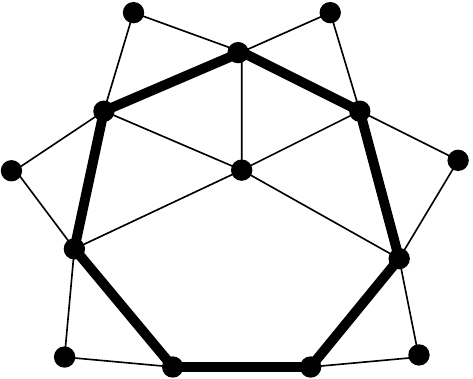} \ &\
\includegraphics[scale=0.35]{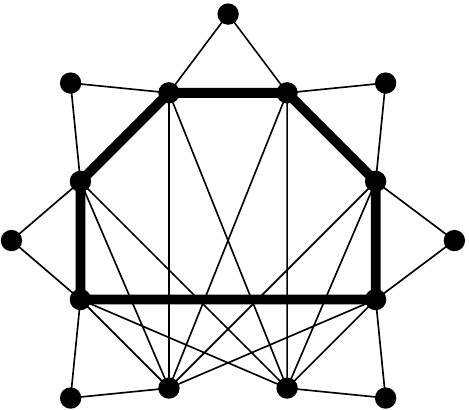} \\
$F_{11}(4k)_{k\geq 2}$ & $F_{12}(4k)_{k\geq 2}$ &
$F_{13}(4k+1)_{k\geq 2}$ &
 $F_{14}(4k+1)_{k\geq 2}$ & $F_{15}(4k+2)_{k\geq 2}$ &
$F_{16}(4k+3)_{k\geq 2}$ \\
\end{tabular}
\caption{Minimal forbidden induced subgraphs for VPT graphs (the
vertices in the cycle marked by bold edges form a clique).}
\label{fig:fiscpath}
\end{figure}

\begin{pf} By Theorem \ref{t:contencion}, $G\in [h+1,2,1]$.
Then, by Theorem \ref{t:cromatic} and since $K$ is a principal
clique of $G$, $\chi(B(G/K))=h+1$.

$(i)$  It is clear that $V(B(G/K))\subseteq V(G)-K$. Suppose there
exists $v\in V(G)-K$ such that $v\notin V(B(G/K))$. Thus, by Claim
\ref{l:branch}, $B(G-v/K)=B(G/K)$. Since $G$ is a minimal non
$[h,2,1]$ graph, $G-v\in [h,2,1]$ and, by Theorem
\ref{t:cromatic}, $B(G-v/K)$ is $h$-colorable. Thus, $B(G/K)$ is
$h$-colorable which contradicts the fact that $K$ is a principal
clique of $G$.

$(ii)$ By item $(i)$ we know that $v\in V(B(G/K))$, then
$|N(v)\cap K|\geq 1$. If $|N(v)\cap K|=1$, $v$ will be an isolated
vertex of $B(G/K)$ and $\chi(B(G/K))=\chi(B(G/K)-v)$. But, by
Claim \ref{l:branch} and Theorem \ref{t:cromatic},
$\chi(B(G/K))=\chi(B(G/K)-v)$ $=\chi(B(G-v/K))=h$, which also
contradicts the fact that $K$ is a principal clique of $G$.

$(iii)$ We know that $\chi(B(G/K))=h+1$. Suppose that $B(G/K)$ is
not $(h+1)$-vertex critical, that is, there is $v\in V(B(G/K))$
such that $\chi(B(G/K)-v)=h+1$. Then, since $v\in V(B(G/K))$, by
Claim \ref{l:branch}, $\chi(B(G-v/K))=\chi(B(G/K)-v)=h+1$, which
contradicts the fact that $G$ is a minimal non $[h,2,1]$ graph.

$(iv)$ We will see that if $N(s_1)\cap K= N(s_2)\cap K$ then
$s_1s_2\notin E(B(G/K))$ and $N_{B(G/K)}(s_1)= N_{B(G/K)}(s_2)$,
which contradicts the fact that $B(G/K)$ is $(h+1)$-vertex
critical. Indeed, if $N(s_1)\cap K= N(s_2)\cap K$ then
$s_1s_2\notin E(B(G/K))$ by definition of branch graph. Moreover,
if $s_3\in N_{B(G/K)}(s_1)$ then there exist $k_1,k_2,k_3\in K$
such that: $(i)$ $k_1s_1\in E(G)$, $k_1s_3\in E(G)$; $(ii)$
$k_2s_1\in E(G)$, $k_2s_3\notin E(G)$; $(iii)$ $k_3s_1\notin
E(G)$, $k_3s_3\in E(G)$. And, since $N(s_1)\cap K= N(s_2)\cap K$,
$k_1s_2\in E(G)$, $k_2s_2\in E(G)$, $k_3s_2\notin E(G)$. In
addition, $s_3s_2\notin E(G)$ because in other case there would be
an induced $4$-cycle $\{s_2,k_2,k_3,s_3\}$ in $G$, which
contradicts the fact that $G\in$ VPT (see Figure
\ref{fig:fiscpath}). Hence, $s_3\in N_{B(G/K)}(s_2)$; we have
proven that $N_{B(G/K)}(s_1)\subseteq N_{B(G/K)}(s_2)$. In a
similar way, it is easy to see that $N_{B(G/K)}(s_2)\subseteq
N_{B(G/K)}(s_1)$.\fin\end{pf}

Theorem \ref{t:split sin dom} shows that all VPT minimal non
$[h,2,1]$ graphs are split without dominated stable vertices.\

To prove this theorem we give the following lemma.

\begin{lem}\label{l:clique}  Let $h\geq 3$, let $G$ be a VPT minimal non $[h,2,1]$
and let $K$ be a principal clique of $G$. Then, $K-\{k\}\in
\mathcal{C}(G-k)$, for all $k\in K$.\end{lem}

\begin{pf} Let $\langle\mathcal{P},T\rangle$ be an
$(h+1,2,1)$-representation of $G$ and let $q\in V(T)$ such that
$K=C_q$. We claim that $\langle\mathcal{P},T\rangle$ is full at
$q$. Indeed, suppose, for a contradiction, that
$\langle\mathcal{P},T\rangle$ is not full at $q$. We can assume,
without loss of generality, that if $x$ is an end vertex of a path
$P_v\in \mathcal{P}$ then there exists a path $P_u\in \mathcal{P}$
intersecting $P_v$ only in $x$, in other case the vertex $x$ can
be removed from $P_v$. This implies that any path of $\mathcal{P}$
linking two branches intersects paths contained in those branches.
Hence, since $\langle\mathcal{P},T\rangle$ is not full at $q$,
there exist branches $T_i$ and $T_j$ of $T$ at $q$ which are
linked by no path of $\mathcal{P}$. Then, by Theorem
\ref{t:repre}, we can obtain a new VPT representation
$\langle\mathcal{P'},T'\rangle$ of $G$ with $d_{T'}(q)\leq h$.
Thus, by Lemma \ref{l:grado}, $B(G/C_q)$ is $h$-colorable which
contradicts the fact that $C_q$ is a principal clique of $G$.

Hence, since $\langle\mathcal{P},T\rangle$ is full at $q$, every
pair of branches of $T$ at $q$ are linked by a path of
$\mathcal{P}$. If there exists $k\in C_q$ such that $C_q-\{k\}$ is
not a clique of $G-k$, there must exists $v\in V(G)-C_q$ such that
$v$ is adjacent to all the vertices of $C_q-\{k\}$. Let $T_1$,
$T_2$,.., $T_{h+1}$ be the branches of $T$ at $q$. Assume, without
loss of generality, that $P_k$ links the branches $T_1$ and $T_2$.
Since $v\in V(G)-C_q$, there exists $i$, such that $P_v$ is
contained in $T_i$. And, since $h\geq 3$, there exists a branch
$T_s$, with $s\neq 1,2,i$. Let $P_u$ be the path of $\mathcal{P}$
linking $T_s$ and $T_r$, with $r\neq i$. It is clear that $u\in
C_q$ and $v$ is not adjacent to $u$, which contradicts the fact
that $v$ is adjacent to all the vertices of $C_q-\{k\}$. Thus,
$C_q-\{k\}\in \mathcal{C}(G-k)$.\fin\end{pf}

The following definition will be used in the proof of Theorem
\ref{t:split sin dom}.

\begin{defn} A $\textbf{canonical VPT representation}$ of $G$ is a pair
$\langle \mathcal{P},T \rangle$ where $T$ is a tree whose vertices
are the members of $\mathcal{C}(G)$, $\mathcal{P}$ is the family
$(P_v)_{v\in V(G)}$ with $P_v=\{C\in \mathcal{C}(G): v\in C\}$ and
$P_v$ is a subpath of $T$ for all $v\in V(G)$.\end{defn}

In  \cite{B5} it was proved that every VPT graph admits a
canonical VPT representation.

\begin{thm} \label{t:split sin dom} Let $G$ be a VPT graph and let $h\geq 3$.
If $G$ is a minimal non $[h,2,1]$ graph,  then $G\in$ Split
without dominated stable vertices.\end{thm}

\begin{pf} Case $(1)$: Suppose that $G\in$ Split with split partition
$(S,K)$, and $G$ has dominated stable vertices. Let $\langle
\mathcal{P},T \rangle$ be a canonical VPT representation of $G$,
and let $q\in V(T)$ such that $K=C_q$. Assume that
$N_T(q)=\{q_1,q_2,..,q_k\}$, with $k>h$, and call
$T_1,T_2,...,T_k$ to the branches of $T$ at $q$ containing the
vertices $q_1,q_2,..,q_k$ respectively. It is clear that for each
$q_i$, with $1\leq i\leq k$, there exists $P_{w_i}\in \mathcal{P}$
such that $q_i\in V(P_{w_i})$ and $q\notin V(P_{w_i})$. Notice
that every $w_i\in S$.

Suppose that $S=\{w_1,w_2,..,w_k\}$. Since $G$ has dominated
stable vertices, by item $(iv)$ of Theorem \ref{l:adyacente1} we
can assume, without loss of generality, that $N(w_1)\subsetneqq
N(w_2)$. This means that $w_1$ and $w_2$ are not adjacent in
$B(G/C_q)$; thus, by item $(iii)$ of Theorem \ref{l:adyacente1},
$N_{B(G/C_q)}(w_1)\nsubseteq N_{B(G/C_q)}(w_2)$. Hence, there
exists $l\in V(B(G/C_q))-\{w_1, w_2\}$, such that $l\in
N_{B(G/C_q)}(w_1)-N_{B(G/C_q)}(w_2)$. Since $V(B(G/C_q))=S$ we can
assume that $l=w_3$. Then, by definition of branch graph, there
exists $z\in C_q$ such that $zw_1\in E(G)$, $zw_3\in E(G)$ and,
since $N(w_1)\subsetneqq N(w_2)$, $zw_2\in E(G)$, which implies
that $P_z$ contains the vertices $q_1$, $q_2$ and $q_3$. Then
$P_z$ is not a path. This contradicts the fact that $\langle
\mathcal{P},T \rangle$ is a VPT representation of $G$.\

We conclude that $S'=S-\{w_1,w_2,..,w_k\}\neq \emptyset$. Let
$G'=G-S'$. Notice that $C_q\in \mathcal{C}(G')$ and
$V(B(G'/C_q))=\{w_1,w_2,..,w_k\}$. Since $G$ is a minimal non
$[h,2,1]$ graph, then $G'\in [h,2,1]$ and $\chi(B(G'/C_q))\leq h$.

We claim that there exists an $h$-coloration of $B(G'/C_q)$ such
that if there exists $x\in C_q$ and $w_i,w_j\in
\{w_1,w_2,..,w_k\}$ with $xw_i\in E(G)$, $xw_j\in E(G)$ then $w_i$
and $w_j$ have different colors in $B(G'/C_q)$. $(*)$

Indeed, if $w_i$ and $w_j$ have the same color in $B(G'/C_q)$ then
$w_iw_j\notin E(B(G'/C_q))$. Then we can assume that
$N(w_i)\subseteq N(w_j)$, since, by hypothesis, there exists $x\in
C_q$ such that $xw_i\in E(G)$ and $xw_j\in E(G)$. Which implies
that $w_i$ is an isolated vertex of $B(G'/C_q)$. Therefore, we can
change the color of $w_i$ to either of the $h-1$ remaining colors.
This process can be done as often as necessary until we have the
desired $h$-coloration of $B(G'/C_q)$.

Hence, we consider an $h$-coloration, say $c'$, of $B(G'/C_q)$
satisfying condition $(*)$.\

Now, we give an $h$-coloration, say $c$, of $B(G/C_q)$ as follows:
given $w\in V(B(G/C_q))$, by Lemma \ref{l:branch1}, there exists
$1\leq i\leq k$ such that $P_w$ is contained in $T_i$, we define
$c(w)=c'(w_i)$. Notice that, in particular, $c(w_i)=c'(w_i)$.

We will see that $c$  is a proper coloration of $B(G/C_q)$. That
is, we have to see that if $uv\in E(B(G/C_q))$ then $c(u)\neq
c(v)$. Since $uv\in E(B(G/C_q))$, by Lemma \ref{l:branch1}, $P_u$
and $P_v$ are in different branches of $T$ at $q$ say $T_i$ and
$T_j$. Moreover, there exists $x\in C_q$  such that $xu\in E(G)$
and $xv\in E(G)$, but this implies that $xw_i\in E(G)$ and
$xw_j\in E(G)$. Hence, since our coloration satisfies condition
$(*)$, $c'(w_i)\neq c'(w_j)$. Thus, $c(u)\neq c(v)$. Therefore,
our coloration is proper.\

Thus, we have an $h$-coloration of $B(G/C_q)$ which contradicts
the fact that $C_q$ is a principal clique of $G$. We conclude
that, if $G\in$ Split then $G$ has no dominated stable vertices.\

Case $(2)$: Suppose that $G \notin$ Split. Since $G$ is a minimal
non $[h,2,1]$ graph, by Theorem \ref{t:contencion}, $G\in
[h+1,2,1]$. Let $\langle \mathcal{P},T \rangle$ be an
$(h+1,2,1)$-representation of $G$ and let $q\in V(T)$ such that
$C_q$ is a principal clique of $G$. We know, by item $(i)$ of
Theorem \ref{l:adyacente1}, that $V(B(G/C_q))=V(G)-C_q$. Since
$G\notin$ Split there exist $x,y\in V(B(G/C_q))$ such that $xy\in
E(G)$.

Let $\tilde{G}$ be the graph which has an
$(h+1,2,1)$-representation $\langle \mathcal{P'},T \rangle$, where
$\mathcal{P'}=(P'_v)_{v\in V(G)}$ such that:

$$P'_v =\left\{
\begin{array}{ll}
P_v,  \ \ \ &\mbox{if\ }v\in C_q\\
 \
q_v, \ \ \ &\mbox{if\ }v\in V(G)-C_q, \ \mbox{where\ }q_v \
\mbox{is the vertex of\ } P_v \ \mbox{closest to q\ }.
\end{array}
\right.$$

Notice that $V(\tilde{G})=V(G)$. We claim that $\tilde{G}$ is a
split graph, with split partition $(V(G)-C_q, C_q)$. Indeed, if
$x,y\in V(G)-C_q$ and $xy\in E(\tilde{G})$ then $q_x=q_y$. Thus,
$N_G(x)\cap C_q = N_G(y)\cap C_q$ which contradicts item $(iv)$ of
Theorem \ref{l:adyacente1}. Hence, $\tilde{G}\in$ Split and, by
Lemma \ref{t:principal}, $C_q$ is a principal clique of
$\tilde{G}$.

 On the other hand, we can assume that $N_G(x)\cap C_q
\subsetneqq N_G(y)\cap C_q$, because in other case it would be an
induced $4$-cycle in $G$ which contradicts the fact that $G\in$
VPT (see Figure \ref{fig:fiscpath}). Then, there exists $w\in C_q$
such that $wx\in E(G)$, $wy\in E(G)$. And, since $xy\in E(G)$ then
$P_x$ and $P_y$ are in a same branch of $T$ at $q$. Hence, by the
existence of $w$, $q_y$ lies on the path of $T$ between $q$ and
$q_x$. Which implies that $\tilde{G}$ has dominated stable
vertices. Now it is easy to see that $B(G/C_q)=B(\tilde{G}/C_q)$,
therefore $\tilde{G}\in [h+1,2,1]-[h,2,1]$.

Then, by Case $(1)$, $\tilde{G}$ is not a minimal non $[h,2,1]$
graph. Thus, there exists $v\in V(\tilde{G})$ such that
$(\tilde{G}-v)\in [h+1,2,1]$.\

If $v\in V(B(\tilde{G}/C_q))$, then
$\chi(B(\tilde{G}-v/C_q))=h+1$. Moreover, by Claim \ref{l:branch}
and since $B(\tilde{G}/C_q)=B(G/C_q)$, we have that
$B(\tilde{G}-v/C_q)=B(\tilde{G}/C_q)-v= B(G/C_q)-v=B(G-v/C_q)$.
Hence, $\chi(B(G-v/C_q))=h+1$ which contradicts the fact that $G$
is a minimal non $[h,2,1]$ graph.\

If $v\in C_q$, then, by Lemma \ref{l:clique}, $C_q-v\in
\mathcal{C}(G-v)$; therefore $C_q-v\in \mathcal{C}(\tilde{G}-v)$.
Thus, $\tilde{G}-v\in$ Split with split partition
$(V(G)-C_q,C_q-v)$. Then, by Lemma \ref{t:principal}, $C_q-v$ is a
principal clique of $\tilde{G}-v$. Hence,
$\chi(B(\tilde{G}-v/C_q-v))=h+1$. Moreover, it is easy to see that
$B(\tilde{G}-v/C_q-v)=B(G-v/C_q-v)$; thus
$\chi(B(\tilde{G}-v/C_q-v))=\chi(B(G-v/C_q-v))=h+1$ which
contradicts the fact that $G$ is a minimal non $[h,2,1]$ graph.

We conclude that $G\in$ Split.\fin\end{pf}

In Theorem \ref{l:adyacente1} we give some necessary conditions on
the branch graph with respect to a principal clique of a minimal
non $[h,2,1]$ graph. In Theorem \ref{t:mascondiciones}, using the
fact that all minimal non $[h,2,1]$ graphs are split without
dominated stable vertices and the fact that the central clique of
a split graph is principal, we will give more necessary conditions
for minimal non $[h,2,1]$ graphs.

\begin{thm}\label{t:mascondiciones} Let $G$ be a VPT graph and let $h\geq 3$.
If $G$ is a minimal non $[h,2,1]$ graph with split partition
$(S,K)$ then: $(i)$ for all $k\in K$, $|N(k)\cap S|=2$; $(ii)$
$|E(B(G/K))|=|K|$; $(iii)$ $B(G/K)$ is $(h+1)$-critical.\end{thm}

\begin{pf} By Theorem \ref{t:split sin dom}, $G\in$ Split without
dominated stable vertices. Let $(S,K)$ be a split partition of
$G$. By Lemma \ref{t:principal}, $K$ is a principal clique of $G$.

$(i)$ Since $G\in$ VPT $\cap$ Split without dominated stable
vertices, $|N(k)\cap S|\leq 2$, for all $k\in K$. Suppose there
exists $k\in K$ such that $|N(k)\cap S|< 2$.

By Theorem \ref{t:contencion}, $G\in [h+1,2,1]$. Let $\langle
\mathcal{P},T \rangle$ be an $(h+1,2,1)$-representation of $G$ and
let $q\in V(T)$ such that $K=C_q$. By Lemma \ref{l:clique},
$C_q-\{k\}\in \mathcal{C}(G-k)$

$1.$ If $|N(k)\cap S|=0$: Then $B(G-k/C_q-\{k\})=B(G/C_q)$. Thus,
$\chi(B(G-k/C_q-\{k\}))=\chi(B(G/C_q))=h+1$, which contradicts the
fact that $G$ is a minimal non $[h,2,1]$ graph.\

$2.$ If $|N(k)\cap S|=1$: We will see that $B(G-k/C_q-\{k\})$
$=B(G/C_q)$. It is clear, by item $(ii)$ of Theorem
\ref{l:adyacente1}, that $V(B(G-k/C_q-\{k\}))=V(B(G/C_q))$ and
$E(B(G-k/C_q-\{k\}))\subseteq E(B(G/C_q))$. Let $uv\in
E(B(G/C_q))$ such that $uv\notin E(B(G-k/C_q-\{k\}))$. Since
$|N(k)\cap S)|=1$ we can assume, without loss of generality, that
$\{N(v)\cap C_q\}-\{N(u)\cap C_q\}=\{k\}$. Therefore, since, for
all $k\in C_q$, $|N(k)\cap S|\leq 2$ we have that
$N_{B(G/C_q)}(v)=\{u\}$ then $d_{B(G/C_q)}(v)=1$, which
contradicts the fact that $H$ is $(h+1)$-vertex critical.

$(ii)$ First we will prove that $|E(B(G/K))|\leq |K|$. Let
$e=s_is_j\in E(B(G/K))$. By definition of branch graph, there
exists $k\in K$ such that $ks_i\in E(G)$, $ks_j\in E(G)$. Thus,
for each $e\in E(B(G/K))$ there exists $k\in K$. Hence, by item
$(i)$, $|E(B(G/K))|\leq |K|$. Now we will see that $|K|\leq
|E(B(G/K))|$. Let $k\in K$. By item $(i)$, $|N(k)\cap S|=2$.
Suppose that $N(k)\cap S=\{s_i,s_j\}$, hence $N(s_i)\cap
N(s_j)\neq \emptyset$. Since there are not dominated stable
vertices, then $N(s_i)\nsubseteq N(s_j)$, $N(s_j)\nsubseteq
N(s_i)$. Thus, $s_is_j\in E(B(G/K))$. Hence, for each $k\in K$
there exist $s_i$,$s_j\in S$ such that $s_is_j\in E(B(G/K))$.
Observe that if $\tilde{k}\in K$ such that $\tilde{k}\neq k$, then
$N(\tilde{k})\cap S\neq N(k)\cap S$. Because if $N(\tilde{k})\cap
S= N(k)\cap S$, then $\tilde{k}$ and $k$ are true twins in $G$
which contradicts the fact that $G$ is a minimal non-$[h,2,1]$
graph. Therefore, $|K|\leq |E(B(G/K))|$.

$(iii)$ By item $(iii)$ of Theorem \ref{l:adyacente1}, $B(G/K)$ is
$(h+1)$-vertex critical. Then, $\chi(B(G/K))$ $=h+1$. We want to
see that $B(G/K)$ is $(h+1)$-edge critical, that is,
$\chi(B(G/K)-e)=h$, for all $e\in E(B(G/K))$. By item $(i)$, for
all $k\in K$, $|N(k)\cap S|=2$ then there are not vertices of $K$
of degree $1$. Moreover, $V(B(G/K))=\{s_1,s_2,...,s_n\}$ with
$\{s_1,s_2,...,s_n\}=S$. Let $e=s_is_j$ and let $k\in K$ such that
$ks_i\in E(G)$, $ks_j\in E(G)$. Since there are not dominated
stable vertices, $B(G-k/K-\{k\})=B(G/K)-e$. Then,
$\chi(B(G-k/K-\{k\}))=\chi(B(G/K)-e)=h$, because $G$ is a minimal
non $[h,2,1]$ graph. Hence, $B(G/K)$ is $(h+1)$-edge critical.
Thus, $B(G/K)$ is $(h+1)$-critical.\fin\end{pf}

\section{Building minimal non [h,2,1] graphs}
\label{s:construc}

The construction presented here is similar to that done in
\cite{pia}, and a generalization of that used in \cite{hugo}.
Given a graph $H$ with $V(H)=$$\{v_{1},...,v_{n}\}$, let $G_H$ be
the graph with vertices:

$$V(G_H) =\left\{
\begin{array}{ll}
v_i,  \ \ \ &\mbox{for each\ }1\leq i \leq n\\
 \
v_{ij}, \ \ \ &\mbox{for each\ }1\leq i<j \leq n \ \mbox{such
that\ } v_{i}v_{j}\in
E(H)\\
 \
\tilde{v_i}, \ \ \ &\mbox{for each\ }1\leq i \leq n \ \mbox{with\
} d_H(v_i)=1.
\end{array}
\right.$$

\begin{figure}[h]
\centering{
\includegraphics[height=1.5in,width=2.7in]{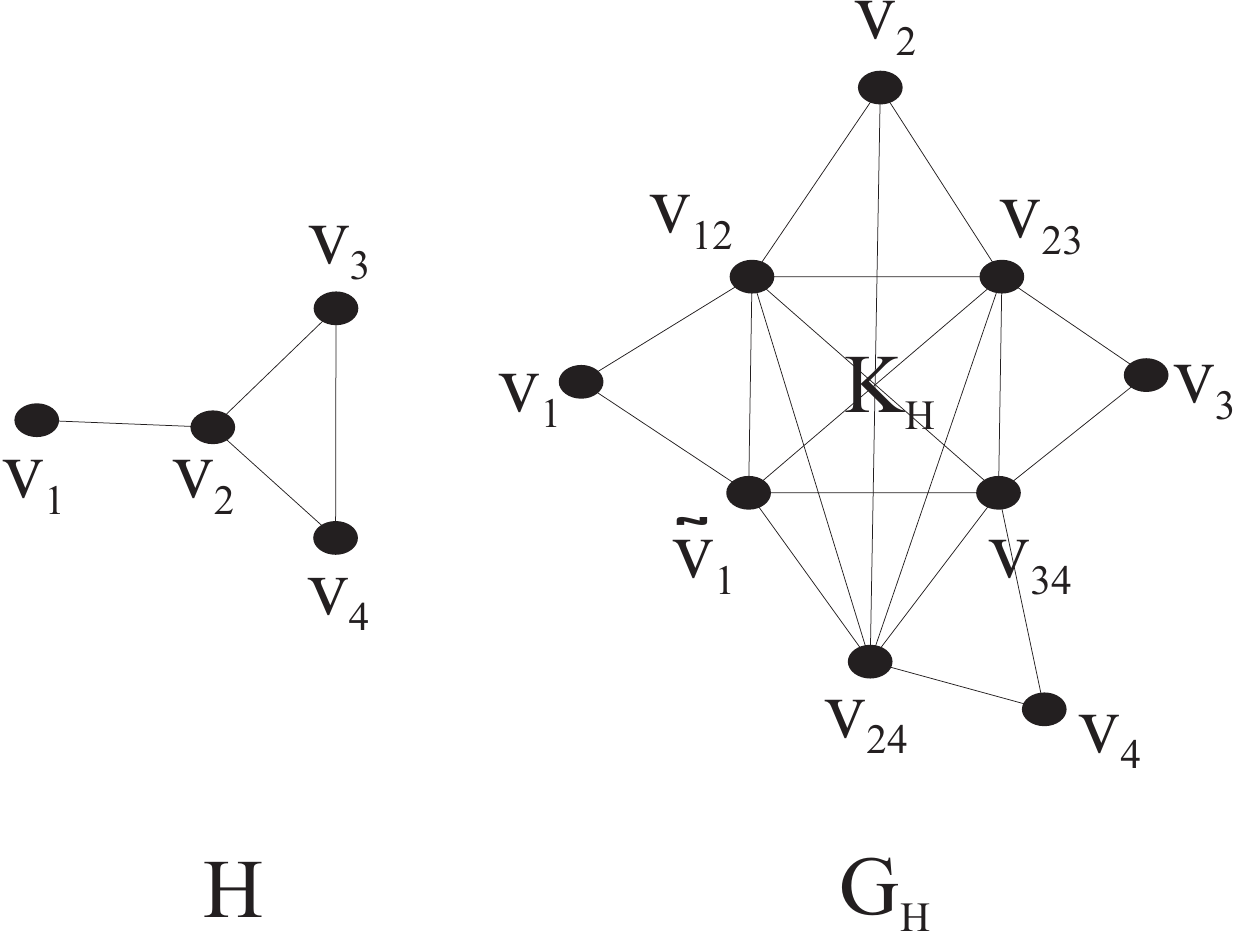}}
\caption{A graph $H$ and the graph $G_{H}$.} \label{f:2}
\end{figure}

The cliques of $G_H$ are: $K_H=\{v_{ij}$, with $1\leq i<j \leq n$
$\}$$\cup$ $\{\tilde{v_i}$, for each $1\leq i \leq n$ such that
$d_H(v_i)=1$$\}$, and $C_{v_i}=\{v_i\}\cup \{v_{ij}; v_j\in
N_H(v_i)\}\cup \{\tilde{v_i},$ if $d_H(v_i)=1\}$, for $1
\leq i \leq n$. 
(See an example in Figure \ref{f:2}).

Notice that the vertices of $G_H$ are partitioned in a stable set
$S_H$ of size $n=|V(H)|$ corresponding to the vertices $v_i$; and
a central clique  $K_H$ of size $|E(H)|+|\{v\in V(H); d_H(v)=1\}|$
corresponding to the remaining vertices. The usefulness of $G_H$
relies on the properties described in the following lemma.

\begin{lem}\cite{pia}\label{l:1-col} $(i)$ $G_H$ is a VPT $\cap$ Split graph without dominated stable
vertices; $(ii)$ $B(G_H/K_H)=H$.\end{lem}

\begin{thm} \label{t:GH minimal} 
Let $h\geq 3$. The graph $G_{H}$ is a minimal non $[h,2,1]$ graph
if and only if $H$ is $(h+1)$-critical.\end{thm}

\begin{pf} Assume that $G_{H}$ is a minimal non $[h,2,1]$ graph.
By item $(ii)$ of Lemma \ref{l:1-col}, $B(G_H/K_H)=H$. Hence, by
item $(iii)$ of Theorem \ref{t:mascondiciones}, $H$ is
$(h+1)$-critical.\

Let $H$ be an $(h+1)$-critical graph with
$V(H)=\{v_1,v_2,...,v_n\}$. By  Lemma \ref{t:principal} and Lemma
\ref{l:1-col}, $Max_{C\in \mathcal{C}(G_{H})}(\chi(B(G_{H}/C)))$
$=\chi(B(G_H/K_H))$ $=\chi(H)=h+1$. Hence, by Theorem
\ref{t:cromatic}, $G_{H}\in [h+1,2,1]-[h,2,1]$. Let us see that
$G_H-v\in [h,2,1]$, for all $v\in V(G_H)$. First, if $v=v_i\in
V(H)$, using Claim \ref{l:branch} and item $(ii)$ of Lemma
\ref{l:1-col}, $B(G_H-v_{i}/K_H)=B(G_H/K_H)-v_{i}=H-v_i$. Thus,
since $H$ is $(h+1)$-vertex critical, $\chi(B(G_H-v_{i}/K_H))=h$.
Hence, $G_H-v_i\in [h,2,1]$. Secondly, if $v=v_{ij}$ being
$e=v_iv_j\in E(H)$, since $B(G_H-v_{ij}/K_H-\{v_{ij}\})=H-e$, then
$\chi(B(G_H-v_{ij}/K_H-\{v_{ij}\}))$ $=\chi(H-e)$. And,
$\chi(H-e)=h$ because $H$ is $(h+1)$-edge critical. Hence,
$G_H-v_{ij}\in [h,2,1]$. Since $H$ has no degree 1 vertices, $G_H$
has no more vertices.\fin\end{pf}

\section{Characterization of minimal non [h,2,1] graphs}

In this Section, we give a characterization of VPT minimal non
$[h,2,1]$ graphs, with $h\geq 3$. The main result of this Section
is Theorem \ref{t:importante} which states that the only VPT
minimal non $[h,2,1]$ graphs are the constructed from
$(h+1)$-critical graphs.

Moreover, in Theorem \ref{t:completo}, we show that the family of
graphs constructed from $(h+1)$-critical graphs together with the
family of minimal forbidden induced subgraphs for VPT
\cite{B9,B11}, is the family of minimal forbidden induced
subgraphs for $[h,2,1]$, with $h\geq 3$.

\begin{thm} \label{t:importante} Let $h\geq 3$ and let $G$ be a VPT graph. $G$ is a minimal non
$[h,2,1]$ graph if and only if there exists an $(h+1)$-critical
graph $H$ such that $G\simeq G_H$.\end{thm}

\begin{pf} The reciprocal implication follows directly applying Theorem \ref{t:GH
minimal}.\

Let $G$ be a minimal non $[h,2,1]$ graph. By Theorem \ref{t:split
sin dom}, we know that $G\in$ Split without dominated stable
vertices. Let $(S,K)$ be a split partition of $G$. By Theorem
\ref{t:contencion}, $G\in [h+1,2,1]$. Let $H=B(G/K)$. By item
$(iii)$ of Theorem \ref{t:mascondiciones}, $H$ is an
$(h+1)$-critical graph. Let us see that $G\simeq G_{H}$. Let
$G_{H}=(S_H,K_H)$. By item $(ii)$ of Lemma \ref{l:1-col},
$B(G_{H}/K_H)=H$ then $B(G_{H}/K_H)=B(G/K)$. Then, since
$V(B(G_H/K_H))=V(B(G/K))$, $S_H=S$. Moreover, since
$E(B(G_H/K_H))=E(B(G/K))$, by item $(ii)$ of Theorem
\ref{t:mascondiciones}, $|K_H|=|K|$ and, by item $(i)$ of Theorem
\ref{t:mascondiciones}, $|N(k)\cap S|=2$ for all $k\in K$. Suppose
that $N(k)\cap S=\{v_i,v_j\}$ we will see that $v_iv_j\in E(H)$.
It is clear that $v_ik\in E(G)$ and $v_jk\in E(G)$. Moreover, by
item $(ii)$ of Theorem \ref{l:adyacente1}, there exist $k',k''\in
K$ such that $k'v_i\in E(G)$, $k''v_j\in E(G)$. If $k'=k''$ then,
since $|N(k)\cap S|=2$ for all $k\in K$, we have that $k'$ and $k$
are true twins in $G$, which contradicts the fact that $G$ is
minimal non $[h,2,1]$ graph. Hence, $k'\neq k''$. Thus,
$k'v_j\notin E(G)$ and $k''v_i\notin E(G)$. Therefore, $v_iv_j\in
E(H)$.\

Hence, we can define a function that assigns to each vertex $k\in
K$ an edge $v_iv_j\in E(H)$, that is, an element of $K_H$. Note
that in $G_H$ the vertex $v_{ij}\in K_H$ is adjacent exactly to
$v_i$ and $v_j$. Hence, the function $f$ can be extended to a new
function $\tilde{f}$ from $K\dot{\cup} S$ to $K_H\dot{\cup} S_H$,
being the identity function from $S$ to $S_H$. Moreover,
$\tilde{f}$ is an isomorphism between $G$ and $G_H$.\fin\end{pf}

\begin{thm}\label{t:completo} Let $h\geq 3$. A graph $G$ is a minimal non $[h,2,1]$ if and only if $G$ is one of the members of
$F_0$, $F_1$,.., $F_{16}$ or $G\simeq G_H$, being $H$ an
$(h+1)$-critical graph.\end{thm}

\begin{pf} By Theorem \ref{t:importante}, if $G\simeq G_H$ being
$H$ an $(h+1)$-critical graph, then $G$ is a minimal non $[h,2,1]$
graph.

If $G$ is any of the members of $F_0$,..,$F_{16}$ then $G\notin$
VPT and $G-v\in$ VPT, for all $v\in V(G)$. Moreover, in
\cite{hugo} it was proved that $G-v\in$ EPT, for all $v\in V(G)$.
Thus, $G-v\in$ VPT $\cap$ EPT$=[3,2,1]$ \cite{B7}, which implies
that $G-v\in [h,2,1]$. Hence, $G$ is a minimal non [h,2,1] graph.

Let $h\geq 3$ and let $G$ be a minimal non $[h,2,1]$ graph.\

Case $(1)$: $G\notin$ VPT. Since $G$ is a minimal non $[h,2,1]$
graph, then $G-v\in [h,2,1]$ for all $v\in V(G)$. Thus, $G-v\in$
VPT for all $v\in V(G)$. Then, $G$ is a minimal forbidden induced
subgraph for VPT. Hence, $G$ is one of the members of $F_0$,
$F_1$,.., $F_{16}$.\

Case $(2)$: $G\in$ VPT. Then, by Theorem \ref{t:importante},
$G\simeq G_H$, being $H$ an $(h+1)$-critical graph.\

Notice that, since every $G_H$ is VPT no member of  $F_0$,
$F_1$,.., $F_{16}$ is an induced subgraph of $G_H$. On the other
hand, suppose that there exists $G$ a member of  $F_0$, $F_1$,..,
$F_{16}$ that has a $G_H$ as induced subgraph. Then, there exists
$\{v_1,..,v_n\}\subseteq V(G)$, with $n\geq 1$, such that
$G-\{v_1,..,v_n\}=G_H$. Hence, since $G_H$ is a minimal non
$[h,2,1]$ graph, $G-\{v_1,..,v_n\}\notin [h,2,1]$ which
contradicts the fact that $G$ is a minimal non $[h,2,1]$
graph.\fin\end{pf}

\begin{thebibliography}{20}\label{bibliography}
\footnotesize{
\bibitem{pia} L. Alc\'on, M. Gutierrez, M. P. Mazzoleni, \emph{Recognizing vertex
intersection graphs of paths on bounded degree trees},
arXiv:1112.3254v1, manuscript. (2011). (Submitted to Discrete
Appl. Math.).

\bibitem{classes} A. Brandst\"{a}dt, V. B. Le, J. P. Spinrad, \emph{Graph Classes: A
Survey}, SIAM Monographs on Discrete Mathematics and Applications.
(1999).

\bibitem{berge} C. Berge, \emph{Graphs and Hypergraphs} (North-Holland, Amsterdam, 1973).

\bibitem{hugo} M. R. Cerioli, H. I. Nobrega, P. Viana, \emph{A
partial characterization by forbidden subgraphs of edge path
graphs}, CTW. (2011) 109-112.

\bibitem{B1} F. Gavril, \emph{The intersection graphs of
subtrees in a tree are exactly the chordal graphs}, J. Combin.
Theory. 16 (1974) 47-56.

\bibitem{gravril} F. Gavril, \emph{A recognition algorithm for the intersection
 graphs of paths in trees}, Discrete Math. 23 (1978) 211-227.

\bibitem{B4} M. C. Golumbic, R. E. Jamison, \emph{Edge and
vertex intersection of paths in a tree}, Discrete Math. 38 (1985)
151-159.

\bibitem{B7} M. C. Golumbic, M. Lipshteyn, M. Stern, \emph{Representing edge
 intersection graphs of paths on degree 4 trees}, Discrete Math. 308 (2008) 1381-1387.

 \bibitem{Jami} R. E. Jamison, H. M. Mulder,
 \emph{Tolerance intersection graphs on binary trees with constant tolerance
3}, Discrete Math. 215 (2000) 115-131.

\bibitem{B6} R. E. Jamison, H. M. Mulder, \emph{Constant tolerance
 intersection graphs of subtress of a tree}, Discrete Math. 290 (2005) 27-46.

\bibitem{B8} C. G. Lekkerkerker, J. C. Boland, \emph{Representation of a
 finite graph by a set of intervals on the real line}, Fund. Math. 51 (1962)
 45-64.

\bibitem{B9} B. L$\acute{e}$v$\hat{e}$que, F. Maffray, M. Preissmann,
\emph{Characterizing path graphs by forbidden induced subgraphs},
J. Graph Theory 62 (2009) 369-384.

\bibitem{B5} C. L. Monma, V. K. Wei, \emph{Intersection graphs
of paths in a tree}, J. Combin. Theory. (1986) 140-181.

\bibitem{nuevo} A. A. Schaffer, \emph{A faster algorithm to recognize undirected path
 graphs}, Discrete Appl. Math. 43 (1993) 261-295.

\bibitem{B11} S. B. Tondato, \emph{Grafos Cordales: Arboles clique y
Representaciones can\'onicas}, Doctoral Thesis, UNLP, Argentina.
(2009).
}
\end{thebibliography}
\end{document}